\documentclass[reqno]{amsart}
\usepackage[latin1]{inputenc}
\usepackage[T1]{fontenc}
\usepackage{amssymb,amsmath,mathrsfs,graphics,color}
\usepackage{hyperref}

\newcommand{\sB}{\mathscr{B}}
\newcommand{\sK}{\mathscr{K}}
\newcommand{\sM}{\mathscr{M}}
\newcommand{\sMin}{\sM^{m}}
\newcommand{\set}[2]{\{#1\mid\,#2\}} 
\newcommand{\R}{\mathbb{R}}
\newcommand{\N}{\mathbb{N}}

\newcommand{\rmax}{\R_{\max}}
\newcommand{\new}[1]{{\em #1}\index{#1}}
\newcommand{\pl}{\flat}
\newcommand{\unit}{0}
\newcommand{\bigunion}{\bigcup}
\newcommand{\union}{\cup}

\newtheorem{prop}{Proposition}[section]

\newtheorem{lemma}[prop]{Lemma}

\newtheorem{theorem}[prop]{Theorem}
\theoremstyle{definition}

\theoremstyle{remark}
\newtheorem{example}[prop]{Example}

\begin{document}

\newcommand\minspace{\sMin}
\newcommand\muminu{\mu^{\rm min}_u}
\newcommand\mumin{\mu^{\rm min}}
\newcommand\mumax{\mu^{\rm max}}
\newcommand\mumaxu{\mu^{\rm max}_u}
\newcommand\sub{m_u}
\newcommand\kay{k}
\newcommand\knu{\kappa^{\nu}}
\newcommand\kmax{\kappa^{\rm max}}

\title{Minimum Representing measures in Idempotent Analysis}
\date{\today}

\begin{abstract}
We show that the set of max-plus measures representing a given max-plus
harmonic vector has a least element. This may be viewed as an
analogue of the uniqueness of the integral representation of harmonic
functions in Potential Theory.
As an application, we show how the distance-like functions of a metric space
can be expressed in terms of the Busemann points of the metric boundary.
\end{abstract}

\maketitle

\section{Introduction}

In Idempotent analysis, one replaces the usual number fields with an idempotent
semifield or semiring. This subject was initiated by Maslov and his collaborators
and developed by many others. An account is given
in~\cite{kolokoltsov_maslov}. Continuing in this tradition,
a recent paper~\cite{AGW-m} develops an idempotent version of potential
theory. The motivation was to find the set of solutions to the
dynamic programming equation of an arbitrary deterministic
Markov decision process with infinite time horizon and infinite
state space. This equation takes the form
\begin{align}
\label{eqn:spectraleqn}
u_i = \sup_{j\in S} (A_{ij} + u_j)
\qquad\mbox{for all $i\in S$,}
\end{align}
where $S$ is the set of states and the map
$A:S\times S\to \R\cup\{-\infty\},\; (i,j)\mapsto A_{ij}$
gives the reward obtained on passing from state $i$ to state $j$.
One searches for solutions 
$u:S\to\R\cup\{-\infty\},\; i\mapsto u_i$.

A fruitful approach to solving~(\ref{eqn:spectraleqn}) is to regard it as
a linear equation in the \emph{max-plus algebra}.
This is the set $\rmax:=\R\cup\{-\infty\}$
equipped with the addition operation $\max(x,y)$
and multiplication operation $x+y$.
Equation~(\ref{eqn:spectraleqn}) is the max-plus analogue of the
equation defining harmonic vectors in Potential Theory,
and for this reason we call its solutions ``max-plus harmonic vectors''.
The analogy is powerful; many results in Potential Theory have analogues.
In particular, one has the following description,
which appeared in~\cite{AGW-m},
of the set of max-plus harmonic vectors.

One first defines the max-plus analogue of the Green kernel:
\begin{equation*}
A^*_{ij}:= \sup\set{A_{ii_1}+\dots+A_{i_{n-1}j}}
            {n\in\N,\; i_1,\dots,i_{n-1}\in S} .
\end{equation*}
This gives the maximal weight of a path from $i$ to $j$.
One assumes that there is some basepoint $b\in S$ from which
every state is accessible, that is,
$A^*_{bj}>-\infty$ for all $j$ in $S$.
The \new{max-plus Martin space} $\sM$ is defined to be the closure of
the set of maps $\sK:=\set{A^*_{\cdot j}- A^*_{bj}}{j\in S}$
in the product topology.
The \new{max-plus Martin boundary} is then defined to be $\sM\setminus\sK$.

Just as in potential theory,
the description of the set of harmonic vectors is in terms of a particular
subset $\sMin$ of the max-plus Martin space, called the
\new{max-plus minimal Martin space}.
(The definition is recalled in the next section.)
It was proved in~\cite{AGW-m} that a vector $u$ is
max-plus harmonic if and only if it is of the form
\begin{align}\label{e-martin}
u=\sup_{w\in \sMin}\mu(w)+ w,
\end{align}
\label{poisson-martin}
where $\mu$ is an upper semicontinuous map from the minimal
max-plus Martin space $\sMin$ to $\rmax$, bounded above.
One may view $\mu$ as
the  max-plus analogue of the density of the \new{spectral measure}
appearing in Potential Theory. The reader may consult~\cite{akian:densities}
for background on max-plus measures.
Unlike their classical counterparts, max-plus measures always
have a density, which can be any upper semicontinuous function.
In this paper, therefore, we will not distinguish between max-plus measures
and their densities.
We introduce a piece of terminology: we say that a max-plus measure
$\nu$ on some subset $D$ of $\sM$ \emph{represents} a vector $u$
if $u=\sup_{\xi\in D}\xi+\nu(\xi)$.

A major difference between the above result and its probabilistic version
is that the representing max-plus measure might not be unique.
A similar degeneracy causes problems when one tries to find representations
for (max-plus) \emph{superharmonic vectors}.
Recall that a vector
$u:S\to\rmax,\; i\mapsto u_i$
is superharmonic if
\begin{align}
\label{eqn:superharmonic}
u_i \ge \sup_{j\in S} (A_{ij} + u_j)
\qquad\mbox{for all $i\in S$.}
\end{align}
If one was looking for a subset $D$ of $\sM$ with the property that every
superharmonic vector has a representing measure defined on $D$, then
a natural choice would be $D:=\sK\union\minspace$, for two reasons.
Firstly, this is an exact analogue of what one does in
(probabilistic) Potential Theory~\cite{dynkin}.
Secondly, as was shown in~\cite{AGW-m}, the elements of
$\sK\union\minspace$ are the normalised extremal vectors,
in the max-plus sense, of the set of superharmonic vectors,
just as in the probabilistic case.
The problem is that it is trivial to find a representing measure on
$\sK\union\minspace$ for every superharmonic vector~$u$: one can take
the upper semicontinuous hull of the map $\sK\to\rmax,\; i\mapsto u_i$.

This problem stems from the fact that there are too many max-plus measures
representing~ $u$. The solution is to be more demanding.
We prove the following theorem.
\begin{theorem}
\label{thm:main}
Assume that $S$ is countable and contains a basepoint from which
every state is accessible.
Let $u\in\rmax^S$ be max-plus superharmonic.
Then, there exists a max-plus measure $\muminu$ on $\sK\union\minspace$
representing~$u$ that is less than any other representing measure
on this set.
If $u$ is max-plus harmonic, then the restriction of $\muminu$ to
$\minspace$ represents $u$ and is less than any other max-plus measure
on $\minspace$ to do so.
\end{theorem}
We call $\muminu$ the \emph{minimum} representing measure of $u$.
One can view the above result as an analogue of the uniqueness of the
spectral measure in Potential Theory.
Note that even if $u$ is max-plus harmonic, $\muminu$ might not take the value
$-\infty$ everywhere on $\sK\setminus\minspace$.
We give an example of this in Section~\ref{sec:example}.

The countability assumption is perhaps not necessary.
However, a proof of the result in the general case would require
different techniques.

To illustrate our ideas, we will show how they apply in an important special
case, that when the kernel $A$ is determined by a metric.
In this setting, the max-plus Martin boundary is just the boundary constructed
by Gromov in~\cite{gromov_hyperbolic}, to which Rieffel has given the
name~\emph{metric boundary}~\cite{rieffel:group}.
This boundary is discussed in~\cite{ballman_gromov_schroeder_manifolds},
\cite{ballmann_spaces}, and~\cite{bridson_haefliger_metric}.
Recent papers concerning it include~\cite{friedland_freitas_pmetrics},
\cite{friedland_freitas_revisiting1},
\cite{karl_metz_nosk_horoballs},
\cite{andreev}, \cite{winweb_boundaries}, and~\cite{winweb_busemann}.

Let $(X,d)$ be a metric space with basepoint. A \emph{distance-like} function
is, according to Gromov~\cite{gromov:hyperbolic}, a function $f$ from $X$
to $\R$ satisfying
\begin{equation}
\label{distancelike}
\inf_{y\in L_t} d(x,y) = f(x) -t
\qquad\text{for all $x\in X$ and $t\le f(x)$,}
\end{equation}
where $L_t := \{x\in X \mid f(x)\le t\}$.
In~\cite{rieffel:group}, Rieffel defines an \emph{almost-geodesic} to be
a function $\gamma$ from some unbounded subset $T$ of $\R_+$ containing $0$
to $X$ such that, for all $\epsilon>0$,
\begin{equation}
|d(\gamma(t),\gamma(s))+d(\gamma(s),\gamma(0))-t| < \epsilon
\end{equation}
for all $t$ and $s$ large enough with $t\ge s$.
He uses the term \emph{Busemann point} to describe (essentially)
the pointwise limit of $d(\cdot,\gamma(t))-d(b,\gamma(t))$ as $t$ tends to
$\infty$, where $\gamma(t)$ is an almost-geodesic.

We use Theorem~\ref{thm:main} to prove the following.
Recall that a proper metric space is one in which the closed balls are compact.
\begin{theorem}
\label{thm:metricmain}
Let $(X,d)$ be a proper geodesic metric space and let $B$ be its set of
Busemann points.
Then, a function $f:X\to \R$ is distance-like if and only if it can be written as
\begin{equation}
\label{infrep}
f=\inf_{h\in B} h +\nu(h),
\end{equation}
where $\nu:B\to \R\union\{\infty\}$ is lower semicontinuous and bounded below.
The set of all such maps $\nu$ satisfying this equation for a fixed distance-like
function $f$ has a greatest element.
\end{theorem}

Section~\ref{sec:preliminaries} recalls some definitions
and Section~\ref{sec:proofs} contains the proofs of the theorems.
We finish with some examples in Section~\ref{sec:example}.

The author would like to thank St\'ephane Gaubert and Marianne Akian for many
stimulating discussions.

\section{Preliminaries}
\label{sec:preliminaries}

We carry over from the Introduction the definitions of the Green kernel
$A^*$, the Martin space $\sM$, and the set $\sK$.
We will also need the following definitions and notation from~\cite{AGW-m}.

The kernel
\begin{equation*}
A^+_{ij}:= \sup\set{A_{ii_1}+\dots+A_{i_{n-1}j}}
            {n\ge 1,\; i_1,\dots,i_{n-1}\in S},
\end{equation*}
gives the maximal weight of a path from $i$ to $j$ of length at least one.
So $A^*_{ij} = A^+_{ij}$ for all $i,j\in S$, $i\neq j$, and $A^*_{ii} = 0$
for all $i\in S$.
We continue to assume that there exists some basepoint $b\in S$ from which
every state is accessible.
The max-plus Martin kernel is
\begin{equation*}
K_{ij}:=A^*_{ij} - A^*_{bj}.
\end{equation*}

For any vector $u$, we use the notation
\begin{align*}
\limsup_{K_{\cdot j} \to w} u_j
   := \inf_{W\ni w} \sup_{K_{\cdot j} \in W} {u}_j,
\end{align*}
where the infimum is taken over all open neighborhoods $W$ of $w$ in $\sM$.
Likewise,
\begin{align*}
\liminf_{K_{\cdot j} \to w} u_j
   := \sup_{W\ni w} \inf_{K_{\cdot j} \in W} {u}_j.
\end{align*}

We say that a sequence $(i_n)_{n\in\N}$ in $S$ converges to $w\in\sM$ if
$K_{\cdot i_n}$ converges pointwise to $w$.

The max-plus minimal Martin space referred to in the Introduction is
defined to be
\begin{align*}
 \sMin:= \set{w\in \sM}{H^{\pl}(w,w)=\unit},
\end{align*}
where
\begin{align*}
H^{\pl}(z,w):=
\limsup_{K_{\cdot i}\to z} \, \liminf_{K_{\cdot j}\to w}
A^*_{bi}+ A^+_{ij} - A^*_{bj}
\qquad
\text{for all $z,w\in \sM$}
\enspace .\
\end{align*}

This notion is closely related to another contained in~\cite{AGW-m},
the notion of \emph{almost-geodesic}.
This is a path $(i_l)_{l\ge 0}$ such that, for some $\beta\in\R$,
\begin{equation}
\label{eqn:def_almost}
A^*_{bi_l} \le \beta + A^*_{bi_0} + A_{i_0i_1} + \dots + A_{i_{l-1}i_l}
\qquad
\text{for all $l\ge0$}.
\end{equation}
Note that this is different from Rieffel's definition of the term, given
earlier. The two notions are compared in~\cite[Section 7]{AGW-m}.
We refer to $\beta$ as the parameter of the almost-geodesic.
It was shown in~\cite{AGW-m} that every almost-geodesic converges
to a point in $\sMin$,
and, conversely, if $\sM$ is first countable, then, for any $w\in\sMin$,
there exists an almost-geodesic converging to $w$.

A definition was also given of an almost-geodesic with respect
to a superharmonic right-vector $u$.
This is a path $(i_l)_{l\ge 0}$ such that, for some $\beta\in\R$,
\begin{equation}
\label{eqn:path}
u_0 \le \beta + A_{i_0i_1} + \dots + A_{i_{l-1}i_l} + u_{i_l}
\qquad
\text{for all $l\ge0$}.
\end{equation}

Let $u\in \rmax^S$ be a vector.
The following max-plus measure $\mumaxu:\sM\to \rmax$ played an
important role in~\cite{AGW-m}:
\begin{align*}
\mumaxu(w) := \limsup_{K_{\cdot j} \to w} A^*_{bj} + {u}_j
\qquad\text{for $w\in \sM$}.
\end{align*}
The map $\mumaxu$ is automatically upper semicontinuous
and bounded above by the constant $u_b$.
In~\cite{AGW-m},
it was shown that $\mumaxu|_{\minspace}$, its restriction to $\minspace$,
represents $u$ and moreover that
it is greater than any other representing measure on this set.

We now give a formula for $\muminu$.
First, define a partial order $\preccurlyeq_u$ on $\sM$
by
\begin{equation*}
z \preccurlyeq_u w
\qquad
\text{if}
\qquad
z_j + \mumaxu(z) \le w_j+ \mumaxu(w)
\qquad
\text{for all $j\in S$}.
\end{equation*}
Next, define $\sub:\sK\cup\minspace \to \rmax$,
\begin{equation*}
\sub(\eta) := 
\begin{cases}
-\infty,
   & \mbox{if there exists $\nu\in\sM\setminus\{\eta\}$
         such that $\eta\preccurlyeq_u \nu$}, \\
\mumaxu(\eta),
   & \mbox{otherwise.}
\end{cases}
\end{equation*}
Finally, take the upper semicontinuous hull:
\begin{equation*}
\muminu(\xi) := \limsup_{\eta\to\xi, \,\, \eta\in\sK\cup\minspace} \sub(\eta)
\qquad\mbox{for all $\xi\in\sK\cup\minspace$}.
\end{equation*}
The proof of Theorem~\ref{thm:main} will consist of showing that $\muminu$
defined in this way has the advertised properties.

\section{Proofs}
\label{sec:proofs}

We begin with some lemmas concerning almost-geodesics.
\begin{lemma}
\label{lem:lemmaA}
Let $(i_n)_{n\in\N}$ be an almost-geodesic with parameter
$\epsilon$ converging to $\xi\in \sK\union \minspace$.
Then
\begin{equation*}
\xi(i_b)
\le \sum_{l=0}^{n-1}A_{i_li_{l+1}} + \xi(i_n) + \epsilon
\qquad\text{for all $n\in\N$}.
\end{equation*}
\end{lemma}
\proof
That $(i_n)_{n\in\N}$ is an almost-geodesic means that
\begin{equation*}
A^*_{bi_n}
\le A^*_{bi_0} + \sum_{l=0}^{n-1}A_{i_li_{l+1}} + \epsilon
\qquad\text{for all $n\in\N$}.
\end{equation*}
We combine this with the triangle inequality
$A^*_{bi_0} + A^*_{i_0i_n} \le A^*_{bi_n}$
and the fact that $A^*_{i_mi_n}$ majorises the sum of the weights along any
path from $i_m$ to $i_n$. We get that
\begin{equation*}
A^*_{i_0i_n} - A^*_{bi_n}
\le \sum_{l=0}^{m-1}A_{i_li_{l+1}} + A^*_{i_mi_n} - A^*_{bi_n} + \epsilon
\qquad\text{for all $n,m\in\N$, $m\le n$}.
\end{equation*}
The result follows on taking the limit as $n\to\infty$.
\qed

\begin{lemma}
[change of basepoint]
\label{lem:lemmaZ}
Let $(i_n)_{n\in\N}$ be an almost-geodesic with parameter $\epsilon$
taking $b$ as basepoint.
Let $j\in S$ be such that $A^*_{ji_0}>-\infty$.
Then $(i_n)_{n\in\N}$ is also an almost-geodesic when $j$ is taken to be
the basepoint, and has parameter
$\epsilon+A^*_{bi_0} - A^*_{bj} - A^*_{ji_0}$.
\end{lemma}
\proof
Take Inequality~(\ref{eqn:def_almost}) and use the fact that
$A^*_{ji_n} \le A^*_{bi_n} - A^*_{bj}$
for all $n\in\N$ to get that
\begin{equation*}
A^*_{ji_n}
 \le A^*_{ji_0} + \sum_{l=0}^{n-1}A_{i_li_{l+1}} + (\epsilon
     + A^*_{bi_0} - A^*_{bj} - A^*_{ji_0})
\qquad\text{for all $n\in\N$.}
\end{equation*}
Note that the bracketed parameter cannot be $+\infty$ since
both $A^*_{bj}$ and $A^*_{ji_0}$ are assumed to be finite.
\qed

\begin{lemma}
\label{lem:lemmaD}
Let $(i_n)_{n\in\N}$ be an almost-geodesic
converging to $\xi\in \sK\union \minspace$.
Let $j\in S$ be such that $A^*_{ji_N}>-\infty$ for some $N\in\N$.
Then
\begin{equation*}
\lim_{n\to\infty} A^*_{ji_n} + \xi(i_n) = \xi(j).
\end{equation*}
Furthermore, if $u$ is a superharmonic vector, then
\begin{equation*}
\lim_{n\to\infty} A^*_{ji_n} + u_{i_n} = \xi(j) + \mumaxu(\xi).
\end{equation*}
\end{lemma}
\proof
By Lemma~7.2 of~\cite{AGW-m},
the truncated path $(i_n)_{n\ge N}$ is an almost-geodesic,
and so by Lemma~\ref{lem:lemmaZ}, it is also an almost-geodesic
when the basepoint is changed to $j$.
Let $\epsilon>0$.
We use Lemma~7.2 of~\cite{AGW-m} again to
deduce that, for $n$ large enough, $(i_l)_{l\ge n}$
is an almost-geodesic with parameter $\epsilon$,
with respect to basepoint $j$.
Therefore
\begin{equation}
\label{eqn:dgeo}
A^*_{ji_m} \le A^*_{ji_n} + A^*_{i_ni_m} + \epsilon
\qquad
\text{for $n$ large enough}.
\end{equation}
Having established this inequality, we restore the status of basepoint to $b$.
Subtracting $A^*_{bi_m}$ from both sides of~(\ref{eqn:dgeo}) and taking the
limit as $m\to\infty$, we find that
\begin{equation*}
\xi(j) \le A^*_{ji_n} + \xi(i_n) + \epsilon
\qquad\text{for $n$ large enough}.
\end{equation*}
Since $\epsilon$ is arbitrary,
it follows that
\begin{equation*}
\liminf_{n\to\infty} A^*_{ji_n} + \xi(i_n) \ge \xi(j).
\end{equation*}
That the quantity $\xi(j)$ is also an upper bound on the
$\limsup$ may be obtained rather easily from the triangle
inequality.
This concludes the proof of the first statement.

We know that $u_{i_n} \ge \xi(i_n)+\mumaxu(\xi)$ from
Lemma~3.6 of~\cite{AGW-m}.
But for any $\delta>0$,
\begin{equation*}
\xi(i_n) \ge A^*_{i_ni_m} - A^*_{bi_m} - \delta
\qquad
\text{for $m$ large enough}.
\end{equation*}
We use these two inequalities and~(\ref{eqn:dgeo})
again to deduce that
\begin{equation*}
A^*_{ji_n} + u_{i_n}
   \ge A^*_{ji_m} - A^*_{bi_m} + \mumaxu(\xi) - \epsilon - \delta.
\end{equation*}
Letting $n$ go to infinity and using the fact that $\epsilon$
and $\delta$ are arbitrary gives us the required lower bound
on the $\liminf$.

The upper bound on the $\limsup$ is easy:
\begin{equation*}
\xi(j) + \mumaxu(\xi)
   = \limsup_{K_{\cdot p}\to \xi} A^*_{jp} +u_p
   \ge \limsup_{n\to\infty}A^*_{ji_n} +u_{i_n}.
\end{equation*}
Thus, the second statement is proved.
\qed

\begin{lemma}
\label{lem:lemmaB}
Assume that $\sM$ is first countable. Let $\xi\in\minspace$ and $j\in S$.
Let $u$ be a superharmonic right vector.
Let $\Delta\ge 0$.
Then
\begin{equation}
\label{eqn:delta}
\Delta\ge u_j - \xi(j)-\mumaxu(\xi)
\end{equation}
if and only if, for all $\epsilon>0$,
there exists an almost-geodesic  with respect to $u$, starting at $j$,
converging to $\xi$, and having parameter $\Delta+\epsilon$.
\end{lemma}
\proof
$(\implies)$
Suppose the inequality above holds.
By Proposition~7.6 of~\cite{AGW-m},
one can find an almost-geodesic $(i_n)_{n\in\N}$ converging to $\xi$.
We see from Lemma~\ref{lem:lemmaZ} that this is also an almost-geodesic 
when the basepoint is taken to $j$ and from Lemma~7.2
of~\cite{AGW-m} that we may take its parameter $\epsilon$ to be as small
as we wish.
So
\begin{equation}
\label{eqn:geo}
A^*_{ji_n}
\le A^*_{ji_0} + \sum_{l=0}^{n-1}A_{i_li_{l+1}} + \epsilon
\qquad\text{for all $n\in\N$}.
\end{equation}
Now take $b$ to be the basepoint again.
By Lemma~\ref{lem:lemmaD},
\begin{equation}
\label{eqn:lemd}
A^*_{ji_n} + u_{i_n} \ge \xi(j) + \mumaxu(\xi) - \epsilon
\qquad\text{for $n$ large enough.}
\end{equation}
Let $(i_l)_{-N\le l\le 0}$ be some finite path starting
at $i_{-N}=j$ and ending at $i_0$, such that
\begin{equation}
\label{eqn:short}
\sum_{l=-N}^{-1}A_{i_li_{l+1}} \ge A^*_{ji_0} - \epsilon.
\end{equation}
Combining Inequalities~(\ref{eqn:delta}),~(\ref{eqn:geo}),
(\ref{eqn:lemd}), and~(\ref{eqn:short}),
we get that
\begin{equation}
u_j \le \Delta+ 3 \epsilon + \sum_{l=-N}^{n-1}A_{i_li_{l+1}} + u_{i_n}
\qquad\text{for $n$ large enough.}
\end{equation}
But
\begin{equation*}
\sum_{l=m}^{n-1}A_{i_li_{l+1}}  + u_{i_n}
\le u_{i_m}
\qquad\text{for all $n$ and $m$ such that $-N\le m\le n$.}
\end{equation*}
So,
\begin{equation*}
u_j \le \Delta+ 3 \epsilon + \sum_{l=-N}^{m-1}A_{i_li_{l+1}} + u_{i_m}
\qquad\text{for all $m\ge -N$,}
\end{equation*}
in other words
$(i_m)_{m\ge-N}$ is an almost-geodesic with parameter $\Delta+ 3 \epsilon$
with respect to $u$.
This proves the first part of the lemma since $\epsilon$
can be chosen arbitrarily.

$(\Longleftarrow)$
Suppose now that, for any $\epsilon> 0$, an almost-geodesic
$(i_n)_{n\in\N}$ exists with the properties specified in
the statement of the lemma.
Then we have
\begin{equation*}
u_{j} \le A^*_{ji_n} + u_{i_n} + \Delta+ \epsilon
\qquad\text{for all $n\in\N$.}
\end{equation*}
To get the desired inequality,
take the limit as $n\to\infty$ using Lemma~\ref{lem:lemmaD},
and use the fact that $\epsilon$ is arbitrary.
\qed

Observe that if $S$ is countable, then $\sM$ is metrisable and hence first
countable.

In the proofs of the next two lemmas, we will use the following notation:
\begin{align*}
\kmax(i,\xi) &:= \xi(i) + \mumaxu(\xi) , \\
\knu(i,\xi) &:= \xi(i) + \nu(\xi)
\qquad \text{and} \\
\kay(i,j) &:= A^*_{ij} + u(j)
\end{align*}
for all $i,j\in S$ and $\xi\in \sK\union\minspace$.
\begin{lemma}
\label{lem:minimality}
Assume that $\sM$ is first countable. Then,
$\sub$ is less than or equal to any max-plus measure on $\sK\union\minspace$
representing $u$.
\end{lemma}
\proof
We know from Theorem~6.1 of~\cite{AGW-m}
that $\mumaxu$ restricted to $\sK\union\minspace$
is greater than any other max-plus measure on $\sK\union\minspace$
representing $u$. Therefore the lemma
will be proved when we show that no max-plus measure $\nu$
on $\sK\union\minspace$
satisfying $\nu\not\ge\sub$ and $\nu\le\mumaxu$ can represent $u$.

For such a $\nu$, there exists some $\xi\in\sK\union\minspace$
such that $\nu(\xi)<\sub(\xi)$. So $\sub(\xi)>-\infty$.
One consequence of this is that
there does not exist $w'\in\sM\setminus\{\xi\}$
such that $w' \succeq \xi$.
Another is that $\sub(\xi) = \mumaxu(\xi)$.

To avoid having to deal separately with the case where $\xi$ is in $\sK$
and the case where it is in $\minspace$, we use the following trick:
we assume that $A_{ii}=0$ for all $i\in S$. We can do this without losing
generality because the Martin kernel $K_{ij}$ does not depend on the diagonal
entries of $A$, and so neither does $\sK$, nor $\sM$, nor whether or not
a given measure represents a given vector. The only effects of setting the
diagonal entries of $A$ to zero are to make all superharmonic vectors
harmonic and to expand $\minspace$ to include $\sK$.
In particular, $\sK\union\minspace$ remains the same.

So now $\xi\in\minspace$.
By Proposition 7.6 of~\cite{AGW-m},
there exists an almost-geodesic converging to $\xi$
with parameter $\epsilon>0$ as small as we wish.
Lemma~\ref{lem:lemmaA} implies that
\begin{equation*}
\knu(i_0,\xi) \le \epsilon + \sum_{m=0}^{n-1}A_{i_m i_{m+1}} + \knu(i_n,\xi)
\qquad
\mbox{for all $n\in\N$.}
\end{equation*}

Since $\nu$ is upper semicontinuous, and $\xi(i_0)$ is continuous in $\xi$,
we see that $\knu(i_0,\cdot)$ is upper semicontinuous. Therefore,
for any $\delta>0$, there exists a set $G\subset\minspace$
containing $\xi$ that is open in  $\minspace$ and such that
\begin{equation*}
\knu(i_0,\eta) < \knu(i_0,\xi) + \delta
\qquad\mbox{for all $\eta\in G$}.
\end{equation*}
Since $\eta$ is superharmonic,
\begin{equation*}
\knu(i_0,\eta) \ge A^*_{i_0 i_n}  + \knu(i_n,\eta)
\qquad
\mbox{for all $\eta\in\sM$ and $n\in\N$.}
\end{equation*}
We deduce from the previous three inequalities that
\begin{equation*}
\knu(i_n,\eta) < \knu(i_n,\xi) + \delta + \epsilon
\qquad\mbox{for all $\eta\in G$ and $n\in\N$}.
\end{equation*}
By choosing $\delta + \epsilon<\sub(\xi)-\nu(\xi)$, we deduce that
\begin{equation}
\label{eqn:maxoverg}
\sup_{\eta\in G}\knu(i_n,\eta) < \xi(i_n)+\sub(\xi)
= \kmax(i_n,\xi)
\qquad\mbox{for all $n\in\N$}.
\end{equation}

Let
\begin{equation*}
L_n:=
\sup_{w\in \minspace\setminus G} \kmax(i_n,w)
\qquad\mbox{for all $n\in\N$}.
\end{equation*}
For any $\delta>0$, there exists a sequence $(\eta_n)_{n\in\N}$ in
$\minspace\setminus G$ such that
$\kmax(i_n,\eta_n) > L_n - \delta$ for all $n\in\N$.
Let $j\in S$.
By the superharmonicity of $\eta_n$,
\begin{equation*}
\kmax(j,\eta_n) \ge A^*_{ji_n}+\kmax(i_n,\eta_n)
\qquad\mbox{for all $n\in\N$}.
\end{equation*}
Also, Lemma~\ref{lem:lemmaD} tells us that, for any $\epsilon>0$,
\begin{equation*}
\kmax(j,\xi)\le A^*_{ji_n}+\kmax(i_n,\xi)+\epsilon
\qquad\mbox{for $n$ large enough}.
\end{equation*}
Putting these inequalities together, we get that
\begin{equation}
\label{eqn:kjetan}
\kmax(j,\eta_n) \ge L_n - \delta +\kmax(j,\xi) - \epsilon - \kmax(i_n,\xi)
\qquad\mbox{for $n$ large enough}.
\end{equation}

The sequence $(\eta_n)_{n\in\N}$ must have at least one limit point $\eta$,
which will necessarily lie in $\sM\setminus G$
and will therefore differ from $\xi$.
Let $(\eta_{n_l})_{l\in\N}$ be a subsequence of $(\eta_n)_{n\in\N}$ 
converging to $\eta$.
Taking the limit infimum of~(\ref{eqn:kjetan}) along $(\eta_{n_l})_{l\in\N}$ and
using the fact that $\kmax(j,\cdot)$ is upper semicontinuous,
and $\delta$ and $\epsilon$ are arbitrary, we conclude that
\begin{equation*}
\kmax(j,\eta) - \kmax(j,\xi)
   > \liminf_{l\to\infty} (L_{n_l}-\kmax(i_{n_l},\xi))
\qquad\mbox{for all $j\in S$}.
\end{equation*}
Since $\eta\nsucceq_u\xi$, the left-hand side is strictly negative
for some $j\in S$, and hence
\begin{equation*}
\liminf_{l\to\infty} (L_{n_l}-\kmax(i_{n_l},\xi)) <0.
\end{equation*}
So by choosing $l$ large enough we can find $n\in\N$ such that
$L_{n}<\kmax(i_n,\xi)$.

Recall that $\nu\le\mumaxu$, which implies that $\knu\le \kmax$.
Applying this and using~(\ref{eqn:maxoverg}), we see that
\begin{equation*}
\sup_{w\in \minspace}
\knu(i_n,w) < \kmax(i_n,\xi).
\end{equation*}
By Lemma~3.6 of~\cite{AGW-m}, $\kmax(i_n,\xi) \le u(i_n)$,
and so we have shown that $\nu$ does not represent $u$.
\qed

\begin{lemma}
\label{lem:representation}
Assume that $S$ is countable.
Let $u\in\rmax^S$ be a harmonic vector.
Then the restriction $\sub|_{\minspace}$ of $\sub$ to $\minspace$
represents $u$.
\end{lemma}
\proof
Let $(i_n)_{n\in\N}$ be a sequence in $S$ that returns to each
element of $S$ infinitely often. Let $I_n := \cup_{m=0}^{n-1}\{i_m\}$
be the set of states visited up to time $n$.
We use the following inductive procedure to define a sequence
$(j_n)_{n\in\N}$ in $S$ and two sequences $(\epsilon_n)_{n\in\N}$ and
$(\delta_n)_{n\in\N}$ of positive reals.

To initialise, we choose arbitrarily 
$j_0$ in $S$, and $\epsilon_0$ and $\delta_0$ greater than zero.

The induction step consists of the following four sub-steps:
\begin{enumerate}

\item[\textit{Step 1.}]
Define the set of states
\begin{equation}
\label{eqn:zdef}
\begin{split}
Z_{n+1} := \big\{
s\in S  \mid
& \kay(i,s)>\kay(i,j_n)-\epsilon_n
\qquad
\text{for all $i\in I_n$},
\quad \text{and} \\
& 
A^+_{j_ns} + u_s >u_{j_n}-\delta_n
\big\},
\end{split}
\end{equation}

\item[\textit{Step 2.}]
\label{choose_j}
Choose $j_{n+1}\in Z_{n+1}$ such that
\begin{equation*}
\kay(i_n,j_{n+1})>\sup_{s\in Z_{n+1}}\kay(i_n,s)-\frac{1}{n}.
\end{equation*}

\item[\textit{Step 3.}]
\label{itm:pathchoice}
Choose a finite path $p^{(n+1)} := p^{(n+1)}_0, \ldots, p^{(n+1)}_{N_{n+1}}$
of some length $N_{n+1}\ge 1$ that starts at
$j_n = p^{(n+1)}_0$ and ends at $j_{n+1} = p^{(n+1)}_{N_{n+1}}$
such that
\begin{equation}
\label{eqn:finitepath}
\sum_{m=0}^{N_{n+1}-1} A_{p^{(n+1)}_m p^{(n+1)}_{m+1}}
   > u_{j_n} - u_{j_{n+1}} - \delta_n.
\end{equation}

\item[\textit{Step 4.}]
Choose $\epsilon_{n+1}$ such that
\begin{equation}
\label{eqn:newep}
0<\epsilon_{n+1} < \epsilon_n +  \kay(i,j_{n+1})-\kay(i,j_n)
\qquad
\text{for all $i\in I_n$.}
\end{equation}
Choose $\delta_{n+1}$ such that
\begin{equation}
\label{eqn:newdel}
0 < \delta_{n+1}
   < \delta_n + u_{j_{n+1}}-u_{j_n}
   + \sum_{m=0}^{N_{n+1}-1} A_{p^{(n+1)}_m p^{(n+1)}_{m+1}}.
\end{equation}
\end{enumerate}

A few remarks necessary to demonstrate that this construction is possible.
Step~2 requires $Z_{n+1}$ to be non-empty. To see that it is, we use
the fact that, since $u$ is harmonic, there exists $r_{n+1}\in S$
such that
\begin{equation}
\label{eqn:exist}
A^+_{j_nr_{n+1}}+u_{r_{n+1}} > u_{j_n} - \min(\delta_n,\epsilon_n).
\end{equation}
That $r_{n+1}$ satisfies the second inequality in the definition
of $Z_{n+1}$ follows immediately. That it satisfies the first follows 
from combining~(\ref{eqn:exist}) with the inequality
$A^*_{ir_{n+1}}\ge A^*_{ij_n} + A^+_{j_nr_{n+1}}$.
So $r_{n+1}$ is in $Z_{n+1}$.

To see that Step 3 is possible, we need only observe that,
since $j_{n+1}$ was chosen to be in $Z_{n+1}$, it must satisfy
\begin{equation*}
A^+_{j_nj_{n+1}} + u_{j_{n+1}} >  u_{j_n} - \delta_n .
\end{equation*}

In Step 4, the choice of $\epsilon_{n+1}$ is possible again
because $j_{n+1}\in Z_{n+1}$,
and so the right-hand-side of~(\ref{eqn:newep}) is positive.
The choice of $\delta_{n+1}$ is possible because~(\ref{eqn:finitepath})
implies that the right-hand-side of~(\ref{eqn:newdel}) is positive.

Denote by $p$ the path obtained by concatenating the paths $p^{(n+1)};n\in\N$.
For all $n\in\N$, let $t_n := \sum_{m=0}^{n-1}N_{m+1}$
be the place at which the finite sequence $p^{(n+1)}$ starts within $p$.
Adding together the inequalities obtained from~(\ref{eqn:newdel})
by varying $n$ from $0$ to $m$,
we get that
\begin{equation*}
0< \delta_{m+1} < \delta_0 + u_{j_{m+1}}
     - u_{j_0} + \sum_{l=0}^{t_{m+1}-1} A_{p_l p_{l+1}}
\qquad
\mbox{for all $m\in\N$}.
\end{equation*}
In other words, $p$ is an almost-geodesic with respect to the right
vector $u$, with parameter $\delta_0$.
Therefore, by Corollary 7.5 of~\cite{AGW-m}
it converges to some $\xi\in \minspace$.
By Lemma~\ref{lem:lemmaB},
\begin{equation}
\label{eqn:costbound}
\delta_0 \ge u_{j_0} - \kmax(j_0,\xi).
\end{equation}

We must show that there is no $\eta\in \sM$ different from $\xi$
such that $\eta\succeq_u \xi$.
So, let $\eta\in \sM$ be such that $\eta\succeq_u \xi$.
Let $s\in S$.
Choose a strictly increasing sequence $(n_q)_{q\in\N}$ in $\N$
such that $i_{n_q}=s$ for all $q\in\N$.

Now we add together the inequalities obtained
from~(\ref{eqn:newep}) by varying $n$
from $n_q$ to $m$, where $m\ge n_q$. We obtain that
\begin{equation*}
0 < 
   \epsilon_{n_q} + \kay(i,j_{m+1}) + \kay(i,j_{n_q}) - \epsilon_{m+1}
\qquad
\mbox{for all $q\in\N$, $i\in I_{n_q}$, and $m\ge n_q$},
\end{equation*}
and that the right-hand-side is increasing in $m$.
So, taking the limit infimum as $m\to\infty$ using Lemma~\ref{lem:lemmaD} yields
\begin{equation*}
\label{eqn:dropbound}
0< \epsilon_{n_q} + \kmax(i,\xi) - \kay(i,j_{n_q})
   - \limsup_{m\to\infty}\epsilon_{m+1}
\qquad
\mbox{for all $q\in\N$ and $i\in I_{n_q}$}.
\end{equation*}
Since $\eta\succeq_u \xi$ and the last term is non-positive,
\begin{equation}
\label{eqn:ibound}
\kmax(i,\eta) > \kay(i,j_{n_q}) - \epsilon_{n_q}
\qquad \mbox{for all $q\in\N$ and $i\in I_{n_q}$}.
\end{equation}

We apply the same procedure to Inequality~(\ref{eqn:newdel}).
Adding the inequalities from $n=n_q$ to $n=m$ yields that
\begin{equation*}
0 < \delta_{n_q} + u_{j_{m+1}} - u_{j_{n_q}}
   + A^*_{j_{n_q}j_{m+1}} -  \delta_{m+1}
\qquad \mbox{for all $q\in\N$ and $m\ge n_q$}
\end{equation*}
and again that the right-hand-side is increasing in $m$.
Taking the same limit as before gives
\begin{equation*}
0< \delta_{n_q} +  \kmax(j_{{n_q}},\xi) - u_{j_{n_q}}
    - \limsup_{m\to\infty}\delta_{m+1}
\qquad \mbox{for all $q\in\N$},
\end{equation*}
and using $\eta\succeq_u \xi$ again gives
\begin{equation}
\label{eqn:jbound}
\kmax(j_{n_q},\eta) > u_{j_{n_q}} - \delta_{n_q}
\qquad \mbox{for all $q\in\N$}.
\end{equation}

It follows from Equations~(\ref{eqn:ibound}) and~(\ref{eqn:jbound}) that,
for any $q\in\N$, any element $\phi$ of $S$ such that $K_{\cdot \phi}$
is close enough to $\eta$ will be in $Z_{n_q+1}$.
So there exists a sequence $(\phi_q)_{q\in\N}$ in $S$ converging to $\eta$
such that $\phi_q\in Z_{n_q+1}$ for all $q\in\N$.
Since $\kmax(s,\cdot)$ is the upper-semicontinuous hull of $\kay(s,\cdot)$,
we may insist that this sequence satisfy additionally
$\lim_{q\to\infty}\kay(s,\phi_q)=\kmax(s,\eta)$.
For this sequence,
\begin{equation*}
\kay(s,j_{n_q+1}) > \kay(s,\phi_q) - \frac{1}{n_q}
\qquad \mbox{for all $q\in\N$}.
\end{equation*}
Taking the limit as $q\to\infty$,
we see that $\kmax(s,\xi) \ge \kmax(s,\eta)$.

Since this holds for all $s\in S$,
we conclude that $\eta$ and $\xi$ are identical.
So we have proved that no element of $\sM$ distinct from $\xi$
is greater than $\xi$ in the ordering $\preceq_u$.
It follows that $\sub(\xi) = \mumaxu(\xi)$.
Combining this with~(\ref{eqn:costbound}), we get that
\begin{align*}
u_{j_0} &\le \xi(j_0) + \sub(\xi) + \delta_0 \\
&\le \sup_{w\in\minspace} w(j_0) + \sub(w) + \delta_0.
\end{align*}
Since $\delta_0$ and $j_0$ may be chosen arbitrarily,
this establishes that
\begin{equation*}
u_\cdot \le \sup_{w\in\minspace} K_{\cdot w} + \sub(w).
\end{equation*}
The opposite inequality follows from the fact that $\mumaxu$
represents $u$ and that $\sub\le\mumaxu$.
\qed

\proof[Proof of Theorem~\ref{thm:main}]
Suppose $u\in\rmax^S$ is (max-plus) harmonic. Lemma~\ref{lem:representation}
shows that $\sub|_{\minspace}$ represents $u$.
Since $\muminu|_{\minspace}$ lies between $\sub|_{\minspace}$ and
$\mumaxu|_{\minspace}$, and both these measures represent $u$,
it must do so also.
If $\nu$ is a
measure on $\minspace$ that represents $u$, then its upper semicontinuous
hull $\bar\nu$ is a measure on $\sK\union\minspace$ that represents $u$.
So then $\sub\le \bar\nu$ by Lemma~\ref{lem:minimality},
and so $\muminu\le \bar\nu$. Taking the restriction to $\minspace$
yields that $\muminu|_{\minspace}\le\bar\nu|_{\minspace}=\nu$.

In order to make Lemma~\ref{lem:representation} applicable in the case
where $u$ is a superharmonic vector, we will employ a trick used in the proof
of Lemma~\ref{lem:minimality}. Recall that setting the diagonal entries
$A_{ii};\,i\in S$ to zero does not change the Martin kernel $K_{ij}$,
nor $\sM$, nor whether a given measure represents a given vector,
but that the new minimal Martin space ${\minspace}'$ is equal to
$\sK\union\minspace$ and now all superharmonic vectors are harmonic.
So we may apply Lemma~\ref{lem:representation} to conclude that
\begin{equation*}
u = \sup_{w\in{\minspace}'}\sub(w)+w
   = \sup_{w\in\sK\union\minspace}\sub(w)+w,
\end{equation*}
in other words, $\sub$ represents $u$. So $\muminu$ also represents $u$.
It was proved in Lemma~\ref{lem:minimality}
that $\sub$ is less than or equal to any measure on $\sK\union\minspace$
representing $u$. It follows that the same is true of $\muminu$.
\qed

\begin{lemma}
\label{lem:metricuniqueness}
Let $(X,d)$ be a separable metric space with a basepoint and let $u$ be a
max-plus harmonic function with respect to the kernel $A$ defined by
$A_{xy}:= -d(x,y)$ for all $x$ and $y$ in $X$.
Then, there exists a max-plus measure $\muminu$ on $\minspace_A$
representing $u$ that is less than any other such measure.
\end{lemma}
\proof
Unfortunately,
we cannot apply Theorem~\ref{thm:main} directly since $X$ is not a countable
set. To get around this problem, we take a countable dense subset $S$ of $X$,
the existence of which is guaranteed by the assumption that $X$ is separable.
We may assume that the basepoint is in $S$.

Since $u$ is max-plus harmonic, it is $1$-Lipschitz.
So, for all $x\in S$,
\begin{equation*}
u(x) = \sup_{y\in X} A_{xy} + u(y) = \sup_{y\in S} A_{xy} + u(y).
\end{equation*}
Thus $u|_S$ is max-plus harmonic with respect to the kernel
$C:=A|_{S\times S}$.

We now need to investigate the relationship between $\minspace_A$ and
$\minspace_C$, the max-plus minimal Martin spaces associated to,
respectively, $A$ and $C$.
Clearly, any almost-geodesic of the kernel $C$ is also
an almost-geodesic of $A$, and its limit in
$\minspace_C$ is just the restriction to $S$ of its limit in $\minspace_A$.
On the other hand, let $(x_n)_{n\in\N}$ be an almost-geodesic of
$A$ converging to $\xi\in \minspace_A$ and let $(\epsilon_n)_{n\in\N}$
be a sequence of positive real numbers satisfying
$\sum_{n=0}^\infty\epsilon_n<\infty$.
One can find a sequence $(y_n)_{n\in\N}$ in $S$ satisfying
$d(x_n,y_n)<\epsilon_n$ for all $n\in\N$.
It follows that $(y_n)_{n\in\N}$ is an almost-geodesic of $C$
and converges to $\xi|_S$.

So we see that the elements of $\minspace_C$ are exactly the restrictions
to $S$ of the elements of $\minspace_A$. Taking this restriction is a bijection
since all the elements of $\minspace_A$ are Lipschitz and $S$ is dense.
Therefore, for any max-plus measure $\mu$ on $\minspace_C$, the map
$\overline\mu:\minspace_A\to\rmax$, $\xi\mapsto \mu(\xi|_S)$ is a max-plus
measure on $\minspace_A$.

If $\mu$ represents $u|_S$, then we deduce from all of the above that,
for $x\in S$,
\begin{equation*}
u(x) = \sup_{\xi\in\minspace_A} \overline\mu(\xi) + \xi(x).
\end{equation*}
But both sides are Lipschitz in $x$, and so the same equation holds for all
$x$ in the closure of $S$, that is $X$.
Therefore, $\overline\mu$ represents $u$.

The converse is clear: if $\overline\mu$ is a max-plus measure on $\minspace_A$
representing $u$, then $\mu$ represents $u|_S$.

The result now follows immediately from the existence,
given by Theorem~\ref{thm:main},
of a minimum representing measure on $\minspace_B$ for $u|_S$.
\qed

\proof[Proof of Theorem~\ref{thm:metricmain}.]

Since $X$ is proper, it is also separable, and so
its set of horofunctions is metrisable. Therefore, for any horofunction $f$,
there is an unbounded sequence $(z_n)_{n\in\N}$ in $X$ such that
$f_n(\cdot):=  d(\cdot,z_n)-d(b,z_n)$ converges to $f$.
Let $x\in X$ and $t\in \R$ be such that $t\le f(x)$.

Since $f_n(x)\le d(x,y)+f_n(y)$ for $n\in \N$ and $y\in X$,
we have that
\begin{equation*}
f(x)\le d(x,y)+f(y) \le d(x,y)+ t
\qquad\text{for all $y\in L_t$}.
\end{equation*}

Since $X$ is a geodesic space, we can find, for each $n\in\N$,
a geodesic line segment $\gamma_n:[0,d(x,z_n)]\to X$ from $x$ to $z_n$.
So, when $d(x,z_n)\ge f(x)-t$, we can take a point $y_n$ in $X$ such that
\begin{equation}
\label{eqn:midpoint}
d(x,z_n) - d(y_n,z_n) = d(x,y_n) = f(x)-t.
\end{equation}
All these points lie inside the closed ball of radius $f(x)-t$ centered at $x$,
which is compact by assumption.
So, by taking a subsequence if necessary, we may assume that $(y_n)_{n\in\N}$
converges to some point $y\in X$. The continuity of the distance function
gives that $d(x,y)=f(x)-t$.
It remains to prove that $y\in L_t$. We have that
\begin{equation*}
d(x,z_n) - d(b,z_n) \to f(x)
\qquad\text{and}\qquad
d(y,z_n) - d(b,z_n) \to f(y)
\end{equation*}
as $n$ tends to $\infty$. Combining these with~(\ref{eqn:midpoint}) gives that
$f(y) = t$.

So, we have proved that every horofunction is distance-like.
In particular, this is true for Busemann functions.

Now suppose that $f$ is of the form
\begin{equation*}
f(x)= \inf_{\alpha\in I}f_\alpha(x)
\qquad\text{for all $x\in X$},
\end{equation*}
where $f_\alpha;\,\alpha\in I$ is some family of distance-like functions on $X$.
Then,
\begin{equation*}
f(x)= \inf_{\alpha\in I}\inf_{y\in L_t^\alpha}d(x,y)+t
   = \inf_{y\in \bigunion_{\alpha\in I}L_t^\alpha}d(x,y)+t
\qquad\text{for all $x\in X$},
\end{equation*}
where $L_t^\alpha:=\{y\in X: f_\alpha(y)\le t\}$
are the corresponding level sets.
For all $t\in \R$ and $\epsilon>0$ such that $t+\epsilon< f(x)$, we have that
\begin{equation*}
\bigunion_{\alpha\in I}L_t^\alpha
   \subset L_t
   \subset \bigunion_{\alpha\in I}L_{t+\epsilon}^\alpha,
\end{equation*}
and so
\begin{equation*}
f(x) \ge \inf_{y\in L_t}d(x,y)+t
     \ge \inf_{y\in \bigunion_{\alpha\in I}L_{t+\epsilon}^\alpha}d(x,y)+t
     = f(x)-\epsilon.
\end{equation*}
Since $\epsilon$ can be chosen as small as we like, we see that $f$ is
distance-like.

So we have proved that the set of distance-like functions is closed under
arbitrary infima. It is also closed under addition of a constant.
It follows that every function of the form given in~(\ref{infrep})
is distance-like.

Now assume that a distance-like function $f:X\to\R$ is given.
We will show that $-f$ is max-plus harmonic with respect to the kernel
$A$ defined by $A_{xy}:= -d(x,y)$ for all $x$ and $y$ in $X$.
Let $x$ and $y$ be in $X$. If $f(x)<f(y)$, then
\begin{equation}
\label{eqn:fupper}
-f(x)\ge A_{xy}-f(y)
\end{equation}
since $A$ is non-positive. On the other hand, if $f(x)\ge f(y)$,
then~(\ref{distancelike}) holds with $t=f(y)$. It follows that
$f(x)\le d(x,y) +f(y)$, and so~(\ref{eqn:fupper}) holds in this case also.
So $-f$ is max-plus superharmonic.

That $-f$ is max-plus subharmonic follows immediately from the fact that
$A_{xx}=0$ for all $x\in X$.

Let $\epsilon>0$ and let $(\epsilon_n)_{n\in\N}$ be a sequence of positive
reals such that $\sum_{i=0}^\infty \epsilon_n<\epsilon$.
Starting at any $x_0\in X$, choose a sequence $(x_n)_{n\in\N}$ in $X$
such that
\begin{equation}
\label{eqn:zed}
f(x_{n+1})\le f(x_{n})-1
\qquad\text{and}\qquad
1\le d(x_{n},x_{n+1})< 1+\epsilon
\end{equation}
for all $n\in\N$.
One can see that this is possible by taking, for successive $n\in\N$,
$t=f(x_n)-1$ in~(\ref{distancelike}) and choosing $x_{n+1}$ in $L_t$
at a distance from $x_n$ close to the infimum.

For all $n\in\N$,
\begin{equation*}
\sum_{i=0}^{n-1} d(x_n,x_{n+1}) \le n+ \sum_{i=0}^{n-1} \epsilon_n
    \le n+ \epsilon
    \le f(x_{0}) - f(x_{n}) + \epsilon.
\end{equation*}
It follows that $(x_n)_{n\in\N}$ is an almost-geodesic in the sense
of~(\ref{eqn:path}) of the kernel $A$
with respect to the function $-f$, having parameter $\epsilon$.

So, by Corollary~7.5 of~\cite{AGW-m}, $(x_n)_{n\in\N}$ converges to some
point $\xi$ in $\minspace$.
From~(\ref{eqn:fupper}) and~(\ref{eqn:zed}),
$d(x_0,x_n)\ge f(x_0)-f(x_n)\ge n$, and so $(x_n)_{n\in\N}$ is unbounded.
It then follows from Proposition~7.10 of~\cite{AGW-m} that $-\xi$ is a
Busemann point. (Beware of the different sign convention used there.)

Returning to the fact that $(x_n)_{n\in\N}$ is an almost-geodesic with respect
to $-f$, we use Lemma~\ref{lem:lemmaB} to deduce that
\begin{equation*}
-f(x_0) \le \xi(x) + \mumax_{-f}(\xi) + \epsilon.
\end{equation*}
But $\epsilon$ is arbitrary, and so
\begin{equation*}
-f(x_0) \le \sup_{\eta\in B} (-\eta(x) + \mumax_{-f}(\eta)).
\end{equation*}
Since the opposite inequality holds by Lemma~3.6 of~\cite{AGW-m}
and since $x_0$ is arbitrary,
\begin{equation*}
f = \inf_{\eta\in B} \eta - \mumax_{-f}(\eta).
\end{equation*}
Recall that the map $\mumax_{-f}|_B$ is upper semicontinuous by construction.
We deduce from the fact that $f(b)=\inf_{\eta\in B}-\mumax_{-f}(\eta)$ that 
$-\mumax_{-f}$ is bounded below.
We have thus proved the first statement of the theorem.

Let $-B:=\{-h:h\in B\}$.
Since $S$ is proper, we may apply Corollary~7.11 of~\cite{AGW-m} to get that
$-B=\minspace\setminus\sK$. We may also deduce from the same fact that
$\sK$ is open in $\sM$, see~\cite{rieffel:group}.
It follows that $-B$ is closed in $\minspace$.
So the map $\mu:\minspace\to\rmax$ defined by
\begin{equation*}
\mu(\xi):=\begin{cases}
   \mumax_{-f}(\xi), & \text{if $\xi\in -B$,} \\
   -\infty, & \text{if $\xi\in \sK$}
   \end{cases}
\end{equation*}
is upper semicontinuous and hence a max-plus measure on $\minspace$.
It represents $-f$ since $\mumax_{-f}|_{-B}$ does.
The fact that $S$ is proper implies that $S$ is also separable, so we may
apply Lemma~\ref{lem:metricuniqueness} to conclude that there exists a max-plus
measure $\mumin_{-f}$ on $\minspace$ representing $-f$ that is smaller than
any other such measure. In particular, $\mumin_{-f}$ must take the value
$-\infty$ outside $-B$. So $\mumin_{-f}|_{-B}$ is a max-plus measure on $-B$
representing $-f$ and it is clearly less than any other such measure
on $-B$ representing $-f$. The greatest lower-semicontinuous map $\nu$
satisfying~(\ref{infrep}) is given by $\mu:=-\mumin_{-f}|_{-B}$.
\qed

\section{Examples}
\label{sec:example}

If $S$ is finite, then it is  known~\cite{finite} that $\muminu(\xi)$
is equal to either $\mumaxu(\xi)$ or $-\infty$ for any superharmonic vector
$u$ and any $\xi\in\sM$. This is not necessarily true however when $S$ is
infinite, as the following example shows.

\begin{example}
We take $S:=\N^2$ and
\begin{equation*}
A_{(x,y)(w,z)}:=
   \begin{cases}
   -1 & \text{if $x=w$ and $z=y\pm 1$}, \\
   -1 & \text{if $(y=0$ or $y=1)$ and $x=w\pm 1,$} \\
   -\infty & \text{otherwise.}
   \end{cases}
\end{equation*}
We choose the basepoint to be $(0,0)$.
A simple calculation shows that the max-plus Martin boundary $\sB$
consists of the infinite set of vectors $a^n;\,n\in\N$
\begin{equation*}
a^n:S\to \rmax,
   (x,y)\mapsto -|x-n|+(2\delta_{xn}-1)|y-1| - 2\delta_{xn}\delta_{y0} +n +1
\end{equation*}
together with the vectors
\begin{equation*}
b^0:S\to \rmax,
   (x,y)\mapsto x-y
\qquad\text{and}\qquad
b^1:S\to \rmax,
   (x,y)\mapsto x-|y-1|+1.
\end{equation*}
Here $\delta_{xn}$ is the Kronecker delta, which takes the value $1$ when
$x=n$ and zero otherwise.
All these points are in $\minspace$.

Define the vector
\begin{equation*}
u:S\to \rmax,
   (x,y)\mapsto \begin{cases}
      x-y & \text{if $y=0$ or $y=1$,} \\
      x+y-4 & \text{if $y\ge 2$.}
   \end{cases}
\end{equation*}
It is easy to see that $u$ is max-plus superharmonic.
One calculates that
\begin{align*}
\mumaxu(a^n) &= -4 \qquad \text{for all $n\in\N$}, \\
\mumaxu(b^0) &= 0,  \\
\mumaxu(b^1) &= -2. 
\end{align*}
The vectors $a^n;\,n\in\N$ are all maximal with respect to the partial
ordering $\preceq_u$, and so $\sub$ agrees with $\mumaxu$ on these vectors.
The same is true for $b^0$. On the other hand, $b^1\preceq_u b^0$ and so
 $\sub(b^1)=-\infty$. However $\muminu(b^1)$ is not $-\infty$ because
$b^1$ is the limit of the sequence $(a^n)_{n\in\N}$.
In fact, $\muminu(b^1)=-4$.
The vectors $a^n;\,n\in\N$ and $b^0$ are isolated, which implies that
\begin{equation*}
\muminu(a^n)= -4 \qquad \text{for all $n\in\N$, and} \qquad
\muminu(b^0) = 0.
\end{equation*}
\qed
\end{example}

Our second example shows that, unlike the situation in probabilistic
potential theory, the minimum representing measure of a max-plus
harmonic function may have weight outside $\minspace$.
\begin{example}
We take $S:=\{1,2,\dots\}\union\{\infty\}$ and
\begin{equation*}
A_{ij}:=
   \begin{cases}
   1/i-1/j & \text{if $j\le i$ and $j<\infty$}, \\
   -1 & \text{if $i=1$ and $j=\infty$}, \\
   -\infty & \text{otherwise.}
   \end{cases}
\end{equation*}
Our basepoint is $\infty$.
A simple calculation reveals that
\begin{equation*}
K_{ij}=
   \begin{cases}
   1/i & \text{if $j\le i$}, \\
   1/i -2 & \text{otherwise}
   \end{cases}
\end{equation*}
and also that
\begin{align*}
\sM &= \{K_{\cdot j} \mid j\in S\}
\qquad\text{and} \\
\minspace &= \{K_{\cdot j} \mid 1\le j <\infty\}.
\end{align*}
Although $K_{\cdot \infty}$ is not in $\minspace$, it is the limit of the
sequence $(K_{\cdot j})_{j\in \N}$, which is contained within $\minspace$.

Consider the vector $u$ defined by $u_i:= 0$ for all $i\in S$.
It is easily verified that $u$ is max-plus harmonic.
One calculates that $\mumaxu(K_{\cdot j}) = -1/j$ for all $j\in S$.
So,
\begin{equation*}
\kmax(i,K_{\cdot j}) = 
   \begin{cases}
   1/i -1/j & \text{if $j\le i$}, \\
   1/i -1/j -2 & \text{otherwise}.
   \end{cases}
\end{equation*}
It follows that $\preceq_u$ is the identity relation, in other words,
no two distinct elements $w$ and $z$ of $\sM$ satisfy $w\preceq_u z$.
So, $\muminu = \sub = \mumaxu$.
\qed
\end{example}

\bibliographystyle{hamsplain}
\bibliography{uniqueness}

\providecommand{\bysame}{\leavevmode\hbox to3em{\hrulefill}\thinspace}
\providecommand{\href}[2]{#2}
\begin{thebibliography}{10}

\bibitem{akian:densities}
Marianne Akian, \emph{Densities of idempotent measures and large deviations},
  Trans. Amer. Math. Soc. \textbf{351} (1999), no.~11, 4515--4543.

\bibitem{AGW-m}
Marianne Akian, St\'ephane Gaubert, and Cormac Walsh, \emph{{The max-plus
  Martin boundary}}, 2004, \mbox{arXiv:math.MG/0412408}, Preprint.

\bibitem{ballmann_spaces}
Werner Ballmann, \emph{Lectures on spaces of nonpositive curvature}, DMV
  Seminar, vol.~25, Birkh\"auser Verlag, Basel, 1995, With an appendix by Misha
  Brin.

\bibitem{ballman_gromov_schroeder_manifolds}
Werner Ballmann, Mikhael Gromov, and Viktor Schroeder, \emph{Manifolds of
  nonpositive curvature}, Progress in Mathematics, vol.~61, Birkh\"auser Boston
  Inc., Boston, MA, 1985.

\bibitem{bridson_haefliger_metric}
Martin~R. Bridson and Andr{\'e} Haefliger, \emph{Metric spaces of non-positive
  curvature}, Grundlehren der Mathematischen Wissenschaften [Fundamental
  Principles of Mathematical Sciences], vol. 319, Springer-Verlag, Berlin,
  1999.

\bibitem{dynkin}
E.B. Dynkin, \emph{Boundary theory of {M}arkov processes (the discrete case)},
  Russian Math. Surveys \textbf{24} (1969), no.~7, 1--42.

\bibitem{friedland_freitas_pmetrics}
Shmuel Friedland and Pedro~J. Freitas, \emph{{$p$}-metrics on {${\rm
  GL}(n,\mathbb{C})/\rm U\sb n$} and their {B}usemann compactifications},
  Linear Algebra Appl. \textbf{376} (2004), 1--18.

\bibitem{friedland_freitas_revisiting1}
\bysame, \emph{Revisiting the {S}iegel upper half plane. {I}}, Linear Algebra
  Appl. \textbf{376} (2004), 19--44.

\bibitem{finite}
St\'ephane Gaubert, Personal communication.

\bibitem{gromov_hyperbolic}
M.~Gromov, \emph{Hyperbolic manifolds, groups and actions}, Riemann surfaces
  and related topics: Proceedings of the 1978 Stony Brook Conference (State
  Univ. New York, Stony Brook, N.Y., 1978) (Princeton, N.J.), Ann. of Math.
  Stud., vol.~97, Princeton Univ. Press, 1981, pp.~183--213.

\bibitem{gromov:hyperbolic}
\bysame, \emph{Hyperbolic groups}, Essays in group theory, Math. Sci. Res.
  Inst. Publ., vol.~8, Springer, New York, 1987, pp.~75--263.

\bibitem{karl_metz_nosk_horoballs}
A.~Karlsson, V.~Metz, and G.~Noskov, \emph{{Horoballs in simplices and
  Minkowski spaces}}, 2004, Preprint.

\bibitem{kolokoltsov_maslov}
Vassili~N. Kolokoltsov and Victor~P. Maslov, \emph{Idempotent analysis and its
  applications}, Mathematics and its Applications, vol. 401, Kluwer Academic
  Publishers Group, Dordrecht, 1997, Translation of {\it Idempotent analysis
  and its application in optimal control} (Russian), ``Nauka'' Moscow, 1994
  Translated by V. E.\ Nazaikinskii, With an appendix by Pierre Del Moral.

\bibitem{andreev}
P.D.Andreev, \emph{{Ideal closures of Busemann space and singular Minkowski
  space}}, 2004, \mbox{arXiv:math.GT/0405121}, Preprint.

\bibitem{rieffel:group}
Marc~A. Rieffel, \emph{Group {$C\sp *$}-algebras as compact quantum metric
  spaces}, Doc. Math. \textbf{7} (2002), 605--651 (electronic).

\bibitem{winweb_boundaries}
Corran Webster and Adam Winchester, \emph{{Boundaries of Hyperbolic Metric
  Spaces}}, 2003, \mbox{arXiv:math.MG/0310101}, {To appear Pacific J.~Math.}

\bibitem{winweb_busemann}
\bysame, \emph{{Busemann Points of Infinite Graphs}}, 2003,
  \mbox{arXiv:math.MG/0309291}, To appear Trans.~Amer.~Math.~Soc.

\end{thebibliography}

\end{document}